\documentclass{article}

\usepackage{arxiv}

\usepackage{algorithm}
\usepackage{algorithmic}
\usepackage{amsmath,amssymb,amsfonts,amsthm}
\usepackage{booktabs}
\usepackage{enumitem}
\usepackage{url}
\usepackage{graphicx}
\usepackage{placeins}
\usepackage{todonotes}
\usepackage{pgfplots}
\usepackage{tikz}
\usepackage{xcolor}
\usepackage{hyperref}
\hypersetup{colorlinks=true,linkcolor=blue,citecolor=blue,urlcolor=blue}

\newtheorem{theorem}{Theorem}

\newtheorem{corollary}[theorem]{Corollary}
\newtheorem{lemma}[theorem]{Lemma}
\theoremstyle{definition}
\newtheorem{definition}[theorem]{Definition}

\theoremstyle{remark}

\newcommand{\R}{\mathbb{R}}
\newcommand{\Xstar}{X^\star}
\newcommand{\fstar}{f^\star}
\newcommand{\dist}{\operatorname{dist}}
\newcommand{\proj}{\operatorname{proj}}
\newcommand{\norm}[1]{\left\lVert#1\right\rVert_2}
\newcommand{\Hdist}{d_{\mathrm H}}
\newcommand{\Rsharpbound}{R_{\mathrm{sharp}}}
\newcommand{\Ren}{R_{\delta}^{\mathrm{env}}}

\title{Bounds on the Distance of Subgradient Trajectories to the Optimal Set}

\author{
Nicholas Dwork$^{1}$, Charles Tsao$^{2}$, Jorge Balb\'as$^{3}$, Cody Coleman$^{4}$, John M. Pauly$^{2}$ \\[3pt]
{\small $^{1}$Department of Biomedical Informatics, University of Colorado School of Medicine} \\
{\small $^{2}$Department of Electrical Engineering, Stanford University}\\
{\small $^{3}$Department of Mathematics, California State University of Northridge}
{\small $^{4}$Department of Computer Science, Stanford University}\\
}

\date{}

\begin{document}
\maketitle

\begin{abstract}
A subgradient method with a constant step size may reach a minimizer and subsequently move away from it. Thus, a bound on the best objective value does not characterize the remainder of the trajectory. In this manuscript, we first establish that the trajectory enters a near-optimal sublevel set in finitely many iterations. We then use the Hausdorff distance between this sublevel set and the optimal set to bound the distance of every subsequent iterate from the optimal set. For objectives with weak sharp minima, we derive a tighter bound and construct examples demonstrating that the bound is tight. Numerical experiments confirm these theoretical results and illustrate when the bounds are conservative.
\end{abstract}

\section{Introduction}

Optimization problems commonly have non-differentiable objective functions; these include those involving absolute values, maxima of affine
functions, and $\ell_1$ regularization.  Applications include compressed sensing \cite{candes2008introduction} and total variation image denoising \cite{rudin1992nonlinear}. For convex problems, the subgradient method offers a simple approach that only requires a subgradient of the objective function at each iteration. A constant step size can make substantial progress initially, but it does not generally produce a convergent
sequence of iterates. For example, a subgradient trajectory for $f(x)=|x|$ can repeatedly move between two points on opposite sides
of its minimizer. Thus, even after the trajectory reaches a point whose objective value is close to optimal, subsequent iterates may move farther from the optimal set.

Classical analyses establish that the constant-step-size subgradient method eventually achieves an objective value close to the optimum \cite{polyak1969minimization,boyd2003subgradient, bertsekas2015convex}. Bertsekas \cite{bertsekas2015convex} showed that the distance to the optimal set decreases whenever the objective error exceeds a threshold determined by the step size and the subgradient bound. Kiwiel \cite{kiwiel2004convergence} studied convergence to neighborhoods of the optimal set for approximate and incremental subgradient methods. More recently, Josz and Lai \cite{josz2023lyapunov} characterized Lyapunov stability near local minima for constant-step subgradient iterations
on locally Lipschitz semi-algebraic functions. They also identified conditions for instability at certain local minima \cite{josz2024sufficient}. In a related paper, they established eventual confinement near a connected component of critical points for several first-order methods on coercive tame functions \cite{josz2024global}. These works address stability in a broader, possibly nonconvex setting.  Instead, our work concerns explicit distances to the global optimal set for convex problems.

A second line of research exploits geometric properties of the objective function to establish stronger convergence guarantees.
Burke and Ferris \cite{burke1993weak} developed the theory of weak sharp minima, a condition relating objective error to distance from the optimal set, and characterized this property for several classes of optimization problems. Nedi\'c and
Bertsekas \cite{nedic2010effect} used sharpness assumptions to analyze subgradient methods with bounded deterministic errors, establishing convergence guarantees even for noncompact constraint sets.
More recently, Johnstone and Moulin \cite{johnstone2020faster} established quantitative bounds for constant-step subgradient methods under H\"olderian growth, which includes weak sharp minima as a special case.
Bai et al.\ \cite{bai2022equivalence} established relationships between weak sharp minima and other error-bound conditions for weakly and approximately convex functions.
Stonyakin et al.\ \cite{stonyakin2023subgradient} studied subgradient methods with Polyak-type step sizes under sharpness assumptions, while Charisopoulos and Davis \cite{charisopoulos2024superlinearly} developed a superlinearly convergent method for sharp semismooth problems.


The preceding studies establish convergence and stability properties of subgradient methods under a variety of assumptions. However, the stability results concern behavior near local minima or critical points, while quantitative bounds on the distance to the optimal set, such as those of Johnstone and Moulin \cite{johnstone2020faster}, rely on additional growth conditions. This motivates our investigation of what can be established about the distance of a constant-step subgradient trajectory to the optimal set without prescribing such a growth condition\footnote{An early version of this work was presented at the Southern California Applied Mathematics Symposium hosted by the Society of Industrial and Applied Mathematicians in 2019 \cite{dwork2019adaptive}.}.
We first establish that the trajectory enters a near-optimal sublevel set in finite time, without assuming a growth condition. We then bound the distance of every subsequent iterate to the optimal set using the Hausdorff distance between the sublevel set and the optimal set. When this Hausdorff distance is finite, the result provides a finite confinement radius. By accounting for the objective error in the subgradient distance inequality, we obtain a potentially tighter bound. Under the additional assumption of weak sharp minima, we derive an explicit asymptotic radius that improves that of Johnstone and Moulin \cite{johnstone2020faster}.

\section{Background}
\label{sec:background}
Let $C\subseteq\R^n$ be a nonempty, closed, convex set, and let $f:\R^n\to\R$ be continuous and convex. We consider the problem $\text{minimize}_{x\in C}\, f(x)$.
The optimal value of the objective function, assumed to be bounded below on $C$, is defined as the following infimum
\begin{equation*}
 \fstar=\inf_{x\in C}f(x)>-\infty.
\end{equation*}
The set of minimizers $\Xstar=\{x\in C:f(x)=\fstar\}$ may be empty. When it is nonempty, $\Xstar$ is closed and convex,
and the optimal value is a minimum.

For a nonempty set $A\subseteq\R^n$, the distance from a point $x\in\R^n$ to $A$ is defined as $\dist(x,A)=\inf_{a\in A}\norm{x-a}$.
For two nonempty sets $A,B\subseteq\R^n$, the (possibly infinite) Hausdorff distance between them is defined as
\begin{equation}
    \Hdist(A,B)=\max\left\{ \sup_{a\in A}\dist(a,B),\ \sup_{b\in B}\dist(b,A) \right\}.
    \label{eq:hausdorff}
\end{equation}

A vector $g\in\R^n$ is a subgradient of $f$ at $x$ if $f(y)\geq f(x)+\langle g,y-x\rangle$ for every $y\in\R^n$.
The subdifferential at $x$ is the collection of all such vectors, and is denoted $\partial f(x)$.

For a nonempty closed set $C\subseteq\R^n$, the Euclidean projection of a point $z\in\R^n$ onto $C$ is defined as $\proj_C(z)=\operatorname{argmin}_{y\in C} \,\norm{z-y}$. Thus, $\proj_C(z)$ is the point in $C$ closest to $z$. When $C$ is also convex, this point exists and is unique.

Given $x_0\in C$ and a constant step size $\alpha>0$, the projected subgradient method generates the trajectory
\begin{equation}
 x_{k+1}=\proj_C(x_k-\alpha g_k),
 \label{eq:update}
\end{equation}
where $g_k\in\partial f(x_k)$ for all $k=0,1,\ldots$. The projection ensures that every iterate is an element of $C$. When $C=\mathbb{R}^n$, \eqref{eq:update} is the ordinary subgradient method.
Throughout the analysis we assume that the selected subgradients satisfy $\norm{g_k}\leq G$ for all $k$ and some fixed $G>0$. Note that this assumption only concerns the selected subgradients along the trajectory; it need not hold for every subgradient at every point in $\R^n$.

For $\delta\geq0$, define the near-optimal sublevel set as $S_\delta=\{x\in C:f(x)\leq\fstar+\delta\}$.
Because $\fstar$ is finite, $S_\delta$ is nonempty for every $\delta>0$, regardless of whether the infimum is attained by $f$.

We begin with a standard inequality \cite{boyd2003subgradient,bertsekas2015convex} in a form that allows comparison with any feasible point, rather than only with a minimizer.

\begin{lemma}[Subgradient comparison inequality]
\label{lem:standard}
Under \eqref{eq:update} and since $\norm{g_k}\leq G$ for all $k$, for any $y\in C$ and every $k\geq0$,
\begin{equation*}
 \norm{x_{k+1}-y}^2\leq\norm{x_k-y}^2 -2\alpha\bigl(f(x_k)-f(y)\bigr)+\alpha^2G^2.
\end{equation*}
\end{lemma}
\begin{proof}
Since $y\in C$, we have $\proj_C(y)=y$.
Projection onto a closed convex set is nonexpansive for any $u,v\in\R^n$: $\norm{\proj_C(u)-\proj_C(v)}\leq\norm{u-v}$.
Therefore,
\begin{align*}
  \norm{x_{k+1}-y}^2 &=\norm{\proj_C(x_k-\alpha g_k)-\proj_C(y)}^2 \leq \norm{x_k-\alpha g_k-y}^2\\
    &=\norm{x_k-y}^2 - 2\alpha\langle g_k,x_k-y\rangle + \alpha^2\norm{g_k}^2
    \leq \norm{x_k-y}^2 - 2\alpha\bigl(f(x_k)-f(y)\bigr) + \alpha^2\norm{g_k}^2
\end{align*}
since the subgradient satisfies $\langle g_k,x_k-y\rangle \geq f(x_k)-f(y)$.
Using $\norm{g_k}\leq G$ proves the claim.
\end{proof}

If $\Xstar\neq\varnothing$, the Euclidean projection $x_k^\star=\proj_{\Xstar}(x_k)$ is well defined. Applying Lemma~\ref{lem:standard} with $y=x_k^\star$ yields
\begin{equation}
 d_{k+1}^2\leq d_k^2
 -2\alpha\bigl(f(x_k)-\fstar\bigr)+\alpha^2G^2,
 \label{eq:fundamental}
\end{equation}
where, for a nonempty set $A$, $\dist(x,A)=\inf_{a\in A}\norm{x-a}$ and $d_k=\dist(x_k,\Xstar)$.

\section{Finite-time entrance and Hausdorff confinement}
\label{sec:hausdorff}

In this section, we establish that the optimization trajectory enters a near-optimal sublevel set in finite time, even when the minimum is not attained, and then derive bounds on its subsequent distance to the optimal set when minimizers exist.

\subsection{Finite-time entrance}

\begin{theorem}[Finite-time entrance]
\label{thm:entrance}
Suppose $C\subseteq\R^n$ is nonempty, closed, and convex, and $f:\R^n\to\R$ is continuous, convex, and bounded below on $C$. Consider the projected subgradient method of Eq.~\eqref{eq:update}, with constant step size $\alpha>0$ and selected subgradients satisfying $\norm{g_k}\leq G$ for every $k$.

For any $\delta > \alpha G^2 / 2$, the trajectory enters the near-optimal sublevel set $S_\delta$ in finitely many iterations.
More precisely, choose any $0< \eta <\delta - \alpha G^2 /2$ and any point $y_\eta\in C$ satisfying $f(y_\eta)<\fstar+\eta$.
Then the first entrance index
$\tau=\min\{k\geq0:x_k\in S_\delta\}$ satisfies
\begin{equation}
 \tau\leq
 \left\lceil
 \frac{\norm{x_0-y_\eta}^2}
 {2\alpha(\delta-\eta)-\alpha^2G^2}
 \right\rceil.
 \label{eq:entranceapprox}
\end{equation}
\end{theorem}

\begin{proof}
Since $\fstar$ is the infimum of $f$ over $C$, we can find a point whose objective value is arbitrarily close to $\fstar$, even when a minimizer does not exist.
In particular, for any $\eta>0$, there exists $y_\eta\in C$ such that $f(y_\eta)<\fstar+\eta$.
Choose $\eta$ as specified in the theorem, and fix one such point $y_\eta$.
We will show that the trajectory cannot remain outside $S_\delta$ indefinitely.

Suppose $x_k\notin S_\delta$. Then $f(x_k)>\fstar+\delta$. Since $f(y_\eta)<\fstar+\eta$, we have $f(x_k)-f(y_\eta)>\delta-\eta$.
Applying Lemma~\ref{lem:standard} with $y=y_\eta$ gives
\begin{align*}
  \norm{x_{k+1}-y_\eta}^2 \leq \norm{x_k-y_\eta}^2 -2\alpha(f(x_k)-f(y_\eta)) +\alpha^2G^2
    <\norm{x_k-y_\eta}^2 - 2\alpha(\delta-\eta)+\alpha^2G^2.
\end{align*}

Define $c_\eta=2\alpha(\delta-\eta)-\alpha^2G^2$. Our choice of $\eta$ guarantees $c_\eta>0$. Thus, whenever $x_k\notin S_\delta$,\\ $\norm{x_{k+1}-y_\eta}^2 < \norm{x_k-y_\eta}^2-c_\eta$.
Therefore, every iteration outside $S_\delta$ decreases the squared distance to $y_\eta$ by more than $c_\eta$.

We now use this fact to show that the trajectory must enter $S_\delta$ in finitely many iterations.
If $x_0\in S_\delta$, then the result is immediate. Otherwise, suppose that the trajectory remains outside $S_\delta$.
Then every iteration outside $S_\delta$ decreases the squared distance to $y_\eta$ by more than $c_\eta>0$.
Thus, if the trajectory remains outside $S_\delta$ for $N$ iterations, then $\norm{x_N-y_\eta}^2 <\norm{x_0-y_\eta}^2-Nc_\eta$.
Choosing $N=\left\lceil \norm{x_0-y_\eta}^2 / c_\eta \right\rceil$ yields a contradiction.
Therefore, the trajectory must enter $S_\delta$ by iteration $N$, proving the claim.
\end{proof}

\begin{corollary}[Finite-time entrance when a minimizer exists]
\label{cor:entranceattained}
Under the assumptions of Theorem~\ref{thm:entrance}, suppose additionally that the optimal set $\Xstar$ is nonempty.
Then, for any $\delta>\alpha G^2/2$, the first entrance index into $S_\delta$ satisfies $\tau \leq \left\lceil d_0^2 / (2\alpha\delta-\alpha^2G^2) \right\rceil$.
\end{corollary}

\begin{proof}
Since $\Xstar$ is nonempty, we can measure the distance from each iterate to the optimal set.
Suppose $x_k\notin S_\delta$. By the definition of $S_\delta$, we have $f(x_k)-\fstar>\delta$.
Therefore, \eqref{eq:fundamental} gives $d_{k+1}^2 < d_k^2-2\alpha\delta+\alpha^2G^2$.
Since $\delta>\alpha G^2/2$, whenever $x_k\notin S_\delta$, $d_{k+1}^2<d_k^2-c$ where $c=2\alpha\delta-\alpha^2G^2$.

In words, every iteration outside $S_\delta$ decreases the squared distance to the optimal set by more than $c$.
We now use this fact to show that the trajectory must enter $S_\delta$ in finitely many iterations.

If $x_0\in S_\delta$, the result is immediate. Otherwise, suppose that the trajectory remains outside $S_\delta$ for its first $N$ iterations.
Repeated application of the above yields $d_N^2<d_0^2-Nc$. Choosing $N=\left\lceil d_0^2 / c \right\rceil$ yields a contradiction.

Therefore, the trajectory must enter $S_\delta$ by iteration $N$. If $N=0$, then $d_0=0$, so $x_0$ is already a minimizer and belongs to $S_\delta$.
\end{proof}

By Theorem \ref{thm:entrance}, the assumption that the infimum is attained is not needed to establish eventual near-optimality.
E.g., let $C=[0,\infty)$ and $f(x)=e^{-x}$. Then $\fstar=0$ but $\Xstar=\varnothing$. Since $f'(x)=-e^{-x}$ has magnitude at most one on $C$, the
selected-subgradient bound holds with $G=1$. The projected iteration is $x_{k+1}=x_k+\alpha e^{-x_k}$ for $x_0\geq0$.
The iterates move toward infinity, but Theorem~\ref{thm:entrance} still guarantees entrance into $S_\delta$ for every $\delta>\alpha/2$. There is no optimal point from which a distance bound could be measured.

\subsection{Confinement when minimizers exist}

For the remainder of this manuscript, we assume that the infimum is attained, so $\Xstar\neq\varnothing$. The set $\Xstar$ is then nonempty, closed, and convex.

Define $H_\delta=\Hdist(S_\delta,\Xstar)$.
Because $\Xstar\subseteq S_\delta$, $H_\delta=\sup_{x\in S_\delta}\dist(x,\Xstar) \ \in[0,\infty]$.
Also define
\begin{equation}
 Q_{\delta,\alpha}=\alpha^2G^2 + \sup_{x\in S_\delta} \left\{\dist(x,\Xstar)^2 -2\alpha\bigl(f(x)-\fstar\bigr)\right\}.
 \label{eq:envelope}
\end{equation}

Unlike an estimate based only on $H_\delta$, the quantity
$Q_{\delta,\alpha}$ retains the objective-error term in
Eq.~\eqref{eq:fundamental}. The following lemma shows
how this quantity bounds the distance to the optimal set
after a subgradient step.

\begin{lemma}[One-step distance bound]
\label{lem:envelope}
Under the assumptions of Corollary~\ref{cor:entranceattained}, if $x_k\in S_\delta$, then $d_{k+1}^2\leq Q_{\delta,\alpha}$.
\end{lemma}
\begin{proof}
Recall from \eqref{eq:fundamental} that $d_{k+1}^2 \leq d_k^2 -2\alpha(f(x_k)-\fstar)+\alpha^2G^2$.
Since $x_k\in S_\delta$,\\ $d_k^2-2\alpha(f(x_k)-\fstar) \leq \sup_{x\in S_\delta} \left\{ \dist(x,\Xstar)^2 -2\alpha(f(x)-\fstar) \right\}$.

Thus $d_{k+1}^2 \leq \alpha^2G^2 + \sup_{x\in S_\delta} \left\{ \dist(x,\Xstar)^2 -2\alpha(f(x)-\fstar) \right\} = Q_{\delta,\alpha}$.
\end{proof}

Because $\Xstar\subseteq S_\delta$, the supremum in \eqref{eq:envelope} is at least zero. Thus
\begin{equation}
 \alpha^2G^2\leq Q_{\delta,\alpha}
 \leq H_\delta^2+\alpha^2G^2,
 \label{eq:envelopeupper}
\end{equation}
where the upper bound is interpreted in the extended real sense if $H_\delta=\infty$.

\begin{theorem}[Hausdorff confinement]
\label{thm:hausdorff}
Assume the hypotheses of Corollary~\ref{cor:entranceattained}.
Let $\delta>\alpha G^2/2$, and define $\Ren~=~\max\left\{H_\delta,\sqrt{Q_{\delta,\alpha}}\right\}$.
If $x_j\in S_\delta$ for some $j\geq0$, then $d_k\leq\Ren$ for every $k\geq j$.
Furthermore, Eq.~\eqref{eq:envelopeupper} gives $\Ren\leq\sqrt{H_\delta^2+\alpha^2G^2}$.
This is a simpler, although potentially less precise, bound on the distance to the optimal set.

By Corollary~\ref{cor:entranceattained}, the trajectory enters $S_\delta$ at some finite iteration $\tau$. Therefore, $d_k\leq\Ren$ holds for every $k\geq\tau$.
\end{theorem}
\begin{proof}
If $H_\delta=\infty$, then $\Ren=\infty$, so the claimed bounds hold in the extended real sense. Instead, assume that $H_\delta<\infty$.

We will show that if $d_k\leq\Ren$, then $d_{k+1}\leq\Ren$. To do this, we consider two cases.
First, suppose $x_k\in S_\delta$. By Lemma~\ref{lem:envelope}, $d_{k+1}^2\leq Q_{\delta,\alpha}$.
Since $\Ren=\max\{H_\delta,\sqrt{Q_{\delta,\alpha}}\}$, we have $Q_{\delta,\alpha}\leq(\Ren)^2$.
Therefore, $d_{k+1}^2\leq Q_{\delta,\alpha}\leq(\Ren)^2$.
Thus, whenever $x_k\in S_\delta$, the next iterate satisfies $d_{k+1}\leq\Ren$.

Now suppose $x_k\notin S_\delta$. By the definition of $S_\delta$, we have $f(x_k)-\fstar>\delta$.
Substituting this inequality into \eqref{eq:fundamental} gives $d_{k+1}^2 < d_k^2-2\alpha\delta+\alpha^2G^2$.
Since $\delta>\alpha G^2/2$, we have $2\alpha\delta-\alpha^2G^2>0$. Therefore, $d_{k+1}^2<d_k^2$.
This means that whenever $x_k$ lies outside $S_\delta$, its distance to the optimal set decreases. In particular, if $d_k\leq\Ren$, then $d_{k+1}\leq\Ren$.

We have now shown that, regardless of whether $x_k$ is inside or outside $S_\delta$, an iterate satisfying $d_k\leq\Ren$ is followed by another iterate satisfying $d_{k+1}\leq\Ren$.

Suppose the trajectory enters $S_\delta$ at iteration $j$. Since $x_j\in S_\delta$, the definition of $H_\delta$ gives $d_j\leq H_\delta\leq\Ren$.
By the argument above, $d_{j+1}\leq\Ren$. By induction, $d_k\leq\Ren$ for all $k\geq j$.
Finally, \eqref{eq:envelopeupper} gives $Q_{\delta,\alpha}\leq H_\delta^2+\alpha^2G^2$.
Consequently, $\Ren = \max\{H_\delta,\sqrt{Q_{\delta,\alpha}}\} \leq\sqrt{H_\delta^2+\alpha^2G^2}$.

By Corollary~\ref{cor:entranceattained}, the trajectory enters $S_\delta$ at some finite iteration $\tau$. Taking $j=\tau$ proves that $d_k\leq\Ren$ for every $k\geq\tau$.
\end{proof}

Theorem~\ref{thm:hausdorff} assumes that the optimal set $\Xstar$ is nonempty, but does not assume that the Hausdorff distance $H_\delta$ is finite. Although the existence of a minimizer ensures that $\Xstar$ is nonempty, it does not guarantee that every point in
$S_\delta$ lies within a bounded distance of $\Xstar$. If $H_\delta$ is finite, the theorem provides a finite bound on the distance of subsequent iterates to the optimal set. If $H_\delta$ is infinite, the theorem remains valid but provides no finite distance bound. The following example illustrates this distinction.

Consider the convex objective $f(u,v)=\sqrt{u^2+v^2}-v$ on $\R^2$. Its optimal value is $\fstar=0$, and
its optimal set is $\Xstar=\{(0,v):v\geq0\}$.
Fix $\delta>0$. For each $t\geq0$, consider the point $y_t=\left(\sqrt{2\delta t+\delta^2},t\right)$.
Evaluating the objective function at this point gives $f(y_t) =\sqrt{2\delta t+\delta^2+t^2}-t =\delta$.
Thus, every point $y_t$ belongs to $S_\delta$.
The closest point in $\Xstar$ to $y_t$ is $(0,t)$.
Therefore, $\dist(y_t,\Xstar)=\sqrt{2\delta t+\delta^2}$.
As $t$ increases, this distance increases without bound, even though every $y_t$ belongs to $S_\delta$. Consequently, $H_\delta =\infty$.
Theorem~\ref{thm:hausdorff} remains valid in this example, but its bound is infinite and therefore provides no finite confinement radius.

\subsection{Behavior as the sublevel set approaches the optimum}

Theorem~\ref{thm:hausdorff} provides a finite confinement radius whenever $H_\delta$ is finite.
The following lemma establishes a sufficient condition for this to occur. Specifically, if the optimal set $\Xstar$ is compact, then every
sublevel set $S_\delta$ is also compact, and therefore $H_\delta$ is finite. The lemma also establishes continuity properties that allow us to obtain an upper bound on $\limsup_{k\to\infty}\dist(x_k,\Xstar)$ by letting $\delta\rightarrow\alpha G^2/2$ in the confinement bound of Theorem~\ref{thm:hausdorff}.

\begin{lemma}[Compact sublevel sets]
\label{lem:compact}
Assume $C$ is nonempty, closed, and convex and $f$ is continuous and convex. Suppose, in addition, that the optimal set $\Xstar$ is nonempty and compact. Then $S_\delta$ is compact for every finite $\delta\geq0$.
\end{lemma}

\begin{proof}
Since $C$ is closed and $f$ is continuous, $S_\delta$ is closed. It remains to show that $S_\delta$ is bounded.

Suppose, for contradiction, that $S_\delta$ is unbounded. Fix $x^\star\in\Xstar$. Then there exists a sequence $z_j\in S_\delta$ such that $\norm{z_j-x^\star}\to\infty$. Define $v_j = (z_j-x^\star)/\norm{z_j-x^\star}$. Since each $v_j$ is a unit vector, compactness of the unit sphere implies that some subsequence converges to a unit vector $v$. We restrict attention to this subsequence.

Fix any $t>0$. For all sufficiently large $j$, we have $t\leq\norm{z_j-x^\star}$. Therefore, $x^\star+tv_j =(1-\lambda_j)x^\star+\lambda_j z_j$ where $\lambda_j=t/\norm{z_j-x^\star}\in[0,1]$. Since $x^\star,z_j\in C$ and $C$ is convex, $x^\star+tv_j\in C$.

Furthermore, convexity of $f$ gives
\begin{equation*}
    f(x^\star+tv_j) \leq (1-\lambda_j)f(x^\star) +\lambda_j f(z_j) \leq (1-\lambda_j)\fstar +\lambda_j(\fstar+\delta)
      = \fstar+\frac{t\delta}{\norm{z_j-x^\star}}.
\end{equation*}

Since $v_j\to v$ and $C$ is closed, $x^\star+tv\in C$. By continuity of $f$, $f(x^\star+tv) = \lim_{j\to\infty}f(x^\star+tv_j) \leq \fstar$.
Since $\fstar$ is the minimum of $f$ over $C$, we conclude that $f(x^\star+tv)=\fstar$. Thus, $x^\star+tv\in\Xstar$.

Because this argument holds for every $t>0$ and $\norm{v}=1$, the optimal set contains minimizers arbitrarily far from $x^\star$. This contradicts the compactness of $\Xstar$. Consequently, $S_\delta$ is bounded. Since it is also closed, $S_\delta$ is compact.
\end{proof}

\begin{lemma}[Right-continuity of sublevel maxima]
\label{lem:sublevelmax}
Under the assumptions of Lemma~\ref{lem:compact}, let $\phi:C\to\mathbb R$ be any continuous function.
Define $M_\phi(\delta)=\max_{x\in S_\delta}\phi(x)$. Then $M_\phi$ is right-continuous in $\delta\geq0$.
\end{lemma}

\begin{proof}
By Lemma~\ref{lem:compact}, $S_\delta$ is nonempty and compact for every $\delta\geq0$. Since $\phi$ is continuous, it attains a
maximum on each $S_\delta$. Thus, $M_\phi(\delta)$ is well defined.
Fix $\delta\geq0$. For each integer $j\geq1$, let $\delta_j=\delta+1/j$. Choose a point $x_j\in S_{\delta_j}$ such that $\phi(x_j)=M_\phi(\delta_j)$.
Since $\delta_j\leq\delta+1$, every $x_j$ belongs to $S_{\delta+1}$. This set is compact, so the sequence $\{x_j\}$ has a convergent subsequence.
Denote this subsequence by $\{x_{j_k}\}$ and its limit by $\bar x$.

By the definition of $S_{\delta_{j_k}}$, $f(x_{j_k})\leq\fstar+\delta+1/j_k$.
Since $f$ is continuous and $x_{j_k}\to\bar x$, taking the limit gives $f(\bar x)\leq\fstar+\delta$.
Also, $\bar x\in C$ because $C$ is closed.
Therefore, $\bar x\in S_\delta$.

Observe that the sets $S_{\delta_j}$ become smaller as $j$ increases. Maximizing a function over a smaller set cannot give a larger maximum. Thus, the sequence $M_\phi(\delta_j)$ is nonincreasing.
Furthermore, since $S_\delta\subseteq S_{\delta_j}$, $M_\phi(\delta_j)\geq M_\phi(\delta)$.
Therefore, $M_\phi(\delta_j)$ converges to some number $L\geq M_\phi(\delta)$.
By continuity of $\phi$,
\begin{equation*}
  L = \lim_{k\to\infty}M_\phi(\delta_{j_k}) =\lim_{k\to\infty}\phi(x_{j_k}) =\phi(\bar x).
\end{equation*}
Since $\bar x\in S_\delta$, we also have $\phi(\bar x)\leq M_\phi(\delta)$.
Combining these inequalities gives $L=M_\phi(\delta)$. Hence, $\lim_{j\to\infty} M_\phi\left(\delta+1/j\right) = M_\phi(\delta)$.

Finally, consider any $s$ satisfying $\delta\leq s\leq\delta + 1/j$.
Because $S_\delta\subseteq S_s\subseteq S_{\delta+1/j}$, we have $M_\phi(\delta) \leq M_\phi(s) \leq M_\phi\left(\delta+\frac{1}{j}\right)$.
Since the upper and lower bounds become arbitrarily close as $j$ increases, $M_\phi(s)$ approaches $M_\phi(\delta)$ as $s$ approaches $\delta$ from above.
Therefore, $M_\phi$ is right-continuous at every $\delta\geq0$.
\end{proof}

\begin{corollary}[Continuity of sublevel-set bounds]
    \label{cor:continuity}
    Under the assumptions of Lemma~\ref{lem:compact}, both $H_\delta$ and $Q_{\delta,\alpha}$ are right-continuous in $\delta\geq0$ for every fixed $\alpha>0$. Furthermore, $\lim_{\delta\to0^+}H_\delta=0$.
\end{corollary}

\begin{proof}
Consider the function $\phi_1(x)=\dist(x,\Xstar)$.
The distance to a nonempty set is continuous, so $\phi_1$ is continuous. By the definition of $H_\delta$, $H_\delta=M_{\phi_1}(\delta)$.
Therefore, Lemma~\ref{lem:sublevelmax} implies that $H_\delta$ is right-continuous.

Next, consider the function $\phi_2(x)=\dist(x,\Xstar)^2 - 2\alpha(f(x)-\fstar)$.
Since both the distance function and $f$ are continuous, $\phi_2$ is also continuous. By the definition of $Q_{\delta,\alpha}$, $Q_{\delta,\alpha} = \alpha^2G^2+M_{\phi_2}(\delta)$.
For fixed $\alpha>0$, the term $\alpha^2G^2$ is constant with respect to $\delta$. Thus, Lemma~\ref{lem:sublevelmax} implies that $Q_{\delta,\alpha}$ is right-continuous.

Finally, $S_0=\Xstar$, so $H_0=\max_{x\in\Xstar}\dist(x,\Xstar)=0$.
Since $H_\delta$ is right-continuous at $\delta=0$, we conclude that $\lim_{\delta\to0^+}H_\delta=H_0=0$.
\end{proof}

\begin{corollary}[Asymptotic Hausdorff bound]
\label{cor:limit}
Under the hypotheses of Lemma~\ref{lem:compact}, also assume the projected subgradient iteration \eqref{eq:update} with constant $\alpha>0$ and bounded selected subgradients. Then
\begin{equation}
 \limsup_{k\to\infty} \hspace{0.5em} d_k\leq \max\left\{H_{\alpha G^2/2}, \sqrt{Q_{\alpha G^2/2,\alpha}}\right\}
    \leq\sqrt{H_{\alpha G^2/2}^2+\alpha^2G^2}.
 \label{eq:limit}
\end{equation}
Furthermore, suppose there exist constants $\alpha_0>0$ and $G>0$ such that, for every step size $0<\alpha\leq\alpha_0$, the selected subgradients satisfy $\norm{g_k}\leq G$ for every $k\geq0$. Then both right-hand bounds tend to zero as $\alpha\to0^+$.
\end{corollary}

\begin{proof}
Fix $\alpha>0$. For any $\delta>\alpha G^2/2$,
Theorem~\ref{thm:hausdorff} guarantees that $d_k\leq\max\left\{ H_\delta,\sqrt{Q_{\delta,\alpha}} \right\}$ for all sufficiently large $k$. Therefore, $\limsup_{k\to\infty}d_k \leq\max\left\{ H_\delta,\sqrt{Q_{\delta,\alpha}} \right\}$.
Since this holds for every $\delta>\alpha G^2/2$, we can let $\delta$ approach $\alpha G^2/2$ from above. By the right-continuity established in Corollary~\ref{cor:continuity}, we obtain $\limsup_{k\to\infty}d_k \leq\max\left\{ H_{\alpha G^2/2}, \sqrt{Q_{\alpha G^2/2,\alpha}} \right\}$.
This establishes the first inequality in \eqref{eq:limit}.

Next, \eqref{eq:envelopeupper} gives $Q_{\delta,\alpha} \leq H_\delta^2+\alpha^2G^2$.
Since $H_\delta\geq0$, both $H_\delta$ and $\sqrt{Q_{\delta,\alpha}}$ are bounded above by $\sqrt{H_\delta^2+\alpha^2G^2}$.
This establishes the second inequality in \eqref{eq:limit}.

Note that $\lim_{\alpha\to0^+}\alpha G^2 = 0$. By Corollary~\ref{cor:continuity}, $\lim_{\alpha\to0^+} H_{\alpha G^2/2}=0$. Therefore, the final expression in \eqref{eq:limit} approaches zero.
\end{proof}

Compactness of $\Xstar$ is sufficient but not necessary for a finite Hausdorff distance. For example, the objective $f(x_1,x_2)=|x_1|$ on $\R^2$ has the unbounded optimal set $\Xstar=\{(0,x_2):x_2\in\R\}$, yet $H_\delta=\delta<\infty$ for every finite $\delta$.

Theorem~\ref{thm:hausdorff} establishes that the subgradient trajectory eventually remains within a bounded distance of the optimal set,
with the confinement radius depending on a selected $\delta$. When the optimal set $\Xstar$ is nonempty and compact, Corollary~\ref{cor:limit} provides an asymptotic bound that does not require choosing a value for $\delta$. However, both bounds involve the quantities $H$ and $Q$,
which may be difficult to calculate. In the next section, we consider an additional condition on the objective function that allows us
to obtain a bound in terms of more accessible problem constants.

\section{A sharper bound for weak sharp minima}
\label{sec:sharp}
The bound of Thm.~\ref{thm:hausdorff} can be made more explicit when the objective grows at least linearly with the distance from its optimal set.

\begin{definition}[Weak sharp minimum]
Assume $\Xstar\neq\varnothing$. The objective function has a
weak sharp minimum with modulus $\mu>0$ on $C$ means $f(x)-\fstar\geq\mu\,\dist(x,\Xstar)$ for all $x\in C$.
\end{definition}

Weak sharpness has been studied extensively in convex optimization \cite{burke1993weak, bai2022equivalence}.
It is the exponent-one case of the H\"olderian growth condition studied in \cite{johnstone2020faster}. Recent studies also analyze
subgradient methods under sharpness assumptions using different
step-size rules \cite{stonyakin2023subgradient, charisopoulos2024superlinearly}. Here, we retain the constant-step projected subgradient iteration, which reduces to the ordinary subgradient method when $C=\R^n$.

Assuming a weak sharp minimum, every point $x\in S_\delta$ satisfies $\dist(x,\Xstar)\leq\delta/\mu$. Hence $H_\delta\leq\delta/\mu$.
This provides an explicit bound on the Hausdorff term in Theorem~\ref{thm:hausdorff}.

\begin{theorem}[Sharp-minimum confinement]
\label{thm:sharp}
Assume that $\Xstar\neq\varnothing$ and that $f$ has a weak sharp minimum on $C$ with modulus $\mu>0$. Consider the projected subgradient iteration \eqref{eq:update} with constant step size $\alpha>0$ and selected subgradients satisfying $\norm{g_k}\leq G$ for all $k\geq0$.
Define $\Rsharpbound = \alpha\max\left\{G,G^2/(2\mu)\right\}$.
If $d_j\leq r$ for some $r\geq\Rsharpbound$, then $d_k\leq r$ for every $k\geq j$.
Additionally,
\begin{equation}
    \limsup_{k\to\infty} \hspace{0.5em} \dist(x_k,\Xstar) \leq \Rsharpbound.
    \label{eq:limsupsharp}
\end{equation}
\end{theorem}

\begin{proof}
By weak sharpness, $f(x_k)-\fstar\geq\mu d_k$.
Substituting this inequality into
\eqref{eq:fundamental} gives
\begin{equation}
d_{k+1}^2\leq
d_k^2-2\alpha\mu d_k+\alpha^2G^2.
\label{eq:sharp_recursion}
\end{equation}

We first show that once the distance to the optimal set is at most $r\geq\Rsharpbound$, it remains at most $r$.
Define $q(t)=t^2-2\alpha\mu t+\alpha^2G^2$. Since $q''(t)=2>0$, the function $q$ is strictly convex. Therefore, its maximum on the interval $[0,r]$ occurs at an endpoint.

Because $r\geq\Rsharpbound$, we have $r\geq\alpha G$ and $r\geq \alpha G^2 /(2\mu)$.
The first inequality gives $q(0)=\alpha^2G^2\leq r^2$.
The second implies that $2\alpha\mu r\geq\alpha^2G^2$, so $q(r)=r^2-2\alpha\mu r+\alpha^2G^2 \leq r^2$.
Since the maximum of $q$ on $[0,r]$ occurs at an endpoint, we conclude that $q(t)\leq r^2$ for every $t\in[0,r]$.
Now suppose $d_k\leq r$. Since distances are nonnegative, $d_k\in[0,r]$. Using \eqref{eq:sharp_recursion}, we obtain $d_{k+1}^2\leq q(d_k)\leq r^2$. Thus $d_{k+1}\leq r$.
By induction, if $d_j\leq r$ for some iteration $j$, then $d_k\leq r$ for every $k\geq j$.

We next show that the trajectory eventually reaches a distance arbitrarily close to $\Rsharpbound$.
Fix any $\varepsilon>0$ and let $r=\Rsharpbound+\varepsilon$.
If $d_k>r$, then
\eqref{eq:sharp_recursion} gives $d_{k+1}^2< d_k^2-\left(2\alpha\mu r-\alpha^2G^2\right)$.
Because $r>\Rsharpbound\geq \frac{\alpha G^2}{2\mu}$, the constant $c=2\alpha\mu r-\alpha^2G^2$ is strictly positive.
Therefore, whenever $d_k>r$, the squared distance decreases by at least $c$. Suppose the trajectory never reaches a point satisfying $d_k\leq r$. After $N$ iterations, $d_N^2<d_0^2-Nc$ for every positive integer $N$.
Choosing $N$ sufficiently large makes the right-hand side negative. This is impossible because $d_N^2\geq0$. Hence the trajectory must reach $d_k\leq r$ after finitely many iterations.

We have already shown that once the distance is at most $r$, it remains at most $r$. Therefore, $\limsup_{k\to\infty}d_k \leq\Rsharpbound+\varepsilon$.
Since this holds for every $\varepsilon>0$, we conclude that $\limsup_{k\to\infty}d_k \leq\Rsharpbound$.
\end{proof}

Theorem~\ref{thm:sharp} also implies that the distance to a neighborhood of the optimal set is non-increasing.
To see this, define $\mathcal N_{\mathrm{sharp}} = \left\{x\in\R^n: \dist(x,\Xstar)\leq\Rsharpbound\right\}$.
The distance from $x_k$ to this neighborhood is $\dist(x_k,\mathcal N_{\mathrm{sharp}}) = \max\{d_k-\Rsharpbound,0\}$.
If $d_k>\Rsharpbound$, then \eqref{eq:sharp_recursion} and $\Rsharpbound\geq\alpha G^2/(2\mu)$ imply $d_{k+1}<d_k$. If $d_k\leq\Rsharpbound$, then the first conclusion of Theorem~\ref{thm:sharp} guarantees that $d_{k+1}\leq\Rsharpbound$.
Therefore, $\dist(x_{k+1},\mathcal N_{\mathrm{sharp}}) \leq \dist(x_k,\mathcal N_{\mathrm{sharp}})$ for every $k\geq0$.
Furthermore, \eqref{eq:limsupsharp} implies that this distance approaches zero as $k\to\infty$.
Thus, the trajectory approaches an invariant neighborhood of the optimal set, even though its distance to the optimal set need not decrease at every iteration.

This is consistent with the familiar behavior of constant-step subgradient methods, which approach a neighborhood of the optimal set without necessarily converging to an optimal point \cite{bertsekas2015convex,johnstone2020faster}.
Theorem~\ref{thm:sharp} identifies an explicit radius for this neighborhood, and the examples below demonstrate that this radius cannot be reduced uniformly under the stated assumptions.

\begin{corollary}[Entrance-time bound]
\label{cor:sharpentrance}
Under the assumptions of Theorem~\ref{thm:sharp}, fix $\varepsilon>0$, and let $K_\varepsilon$ be the first iteration at which $d_{K_\varepsilon} \leq \Rsharpbound+\varepsilon$.
Then
\begin{equation}
    K_\varepsilon\leq \left\lceil \frac{d_0^2}{2\alpha\mu(\Rsharpbound+\varepsilon) - \alpha^2G^2} \right\rceil.
    \label{eq:sharpentrance}
\end{equation}
\end{corollary}

\begin{proof}
Let $r=\Rsharpbound+\varepsilon$ and $c=2\alpha\mu r-\alpha^2G^2$. As established in the proof of Theorem~\ref{thm:sharp}, $c>0$ and $d_{k+1}^2<d_k^2-c$ whenever $d_k>r$,

If the trajectory remains outside
this region for its first $N$ iterations,
repeated application of this inequality
gives $d_N^2<d_0^2-Nc$.
For $N=\left\lceil d_0^2/c\right\rceil$, the right-hand side is nonpositive, contradicting $d_N^2\geq0$. Thus, the trajectory must enter the region by iteration $N$. If $N=0$, then $d_0=0$, so the trajectory is already inside the region. Substituting the definition of $c$ proves the claim.
\end{proof}

We now compare the asymptotic bound of Theorem~\ref{thm:sharp} with an existing result for objectives satisfying weak sharpness. Johnstone and Moulin \cite{johnstone2020faster} established distance bounds under a more general H\"olderian growth condition. Specializing their Theorem~1 to weak sharp minima gives $R_{\mathrm{JM}} = \alpha G\sqrt{1+\left(G/(2\mu)\right)^2}$.
The same expression can be obtained from the simpler Hausdorff-based bound in Theorem~\ref{thm:hausdorff}, using $H_\delta\leq\delta/\mu$ and taking the limit as $\delta\to(\alpha G^2/2)^+$.
In contrast, Theorem~\ref{thm:sharp} gives $\Rsharpbound = \alpha G\max\left\{1,G/(2\mu)\right\} < R_{\mathrm{JM}}$.
The relative reduction in the confinement radius is largest when $G=2\mu$, where it equals $1-1/\sqrt{2}\approx29.3\%$. Note that Johnstone and Moulin also establish convergence-rate results that are not addressed here.

\subsection{Examples where the bound is attained}
\label{sec:attainmentExamples}
The two terms of $\Rsharpbound$ are both necessary. The following examples show that the coefficients cannot be decreased uniformly over the stated class of convex functions.

Consider $f(x)=|x|$, whose optimal set is $\{0\}$.
Starting at $x_0=0$, choose the subgradient $g_0=1$.
At a negative point choose $g=-1$ and at zero choose $g=1$.
The resulting trajectory is $0,-\alpha,0,-\alpha,\ldots$.
Thus $\limsup_kd_k=\alpha=\alpha G=\Rsharpbound$.
The trajectory can also avoid the minimizer: starting at
$x_0=\epsilon\alpha$, for $0<\epsilon<1$, gives the cycle
$\epsilon\alpha,-(1-\epsilon)\alpha,\ldots$, whose largest
distance approaches $\alpha$ as $\epsilon\rightarrow0$.

Separately, consider $f(x)=\|x\|_\infty$ on $\R^n$ where $\|x\|_\infty=\max_{1\leq i\leq n}|x_i|$.
The optimal set is $\Xstar=\{0\}$ and $\fstar=0$. Since $|x_i|\leq\|x\|_\infty$ for every component of $x$, we have $\norm{x}^2=\sum_{i=1}^n x_i^2 \leq n\|x\|_\infty^2$.
Therefore, $f(x)-\fstar=\|x\|_\infty \geq\frac{1}{\sqrt n}\norm{x} =\frac{1}{\sqrt n}\dist(x,\Xstar)$.
Thus, $f$ has a weak sharp minimum with modulus $\mu=1/\sqrt n$.
Furthermore, the infinity norm is Lipschitz continuous with constant one in the Euclidean norm. Consequently, every subgradient $g\in\partial f(x)$ satisfies $\norm{g}\leq1$, so we may take $G=1$.

Set $x_0=(\alpha/2,\ldots,\alpha/2)$. 
At each iteration, choose
\[
g_k= \begin{cases}
(1,0,\ldots,0), & x_{k,1}>0,\\
(-1,0,\ldots,0), & x_{k,1}<0.
\end{cases}
\]
These are valid subgradients because the absolute value of the first coordinate equals $\|x_k\|_\infty$ at every iteration. The update changes the first coordinate from $\alpha/2$ to $-\alpha/2$, or vice versa, while leaving all other coordinates unchanged. This choice is valid because every coordinate has magnitude $\alpha/2$.
Consequently, $d_k = \norm{x_k} = \alpha\sqrt n / 2 = \alpha G^2/2\mu)$ for every $k$.
For $n\geq4$, this is equal to $\Rsharpbound$.

Another example is $f(x)=\|x\|_1$ on $\R^2$. Here $G=\sqrt2$ and $\mu=1$. Starting at the origin and choosing the subgradient $(1,1)$ produces the cycle between $(0,0)$ and $(-\alpha,-\alpha)$, so the largest distance is $\alpha\sqrt2=\alpha G=\Rsharpbound$.

These constructions show that, without additional assumptions, $\Rsharpbound$ cannot be reduced.

\section{Experiments}
\label{sec:experiments}
We performed numerical experiments to examine the entrance-time bound of Corollary~\ref{cor:entranceattained} and the confinement bounds of Theorems~\ref{thm:hausdorff} and \ref{thm:sharp}. For problems with weak sharp minima, we also evaluated the distance bound implied by Johnstone and Moulin \cite{johnstone2020faster}. The trajectories are those of the ordinary subgradient method. All minimizers, subgradients, and bounds used below are known analytically. All experiments used a constant step size and deterministic subgradient selection. The default step size was $\alpha=0.2$. Distances are measured using the Euclidean norm. The theoretically predicted values were calculated from the known constants $G$ and $\mu$.

\subsection{Exact confinement examples}
Table~\ref{tab:exact} reports the largest late-iterate distance
for the exactness examples of Section~\ref{sec:attainmentExamples}. We compare the
improved radius $\Rsharpbound$ to the radius
$R_{\mathrm{JM}}$ implied by Theorem~1 of
\cite{johnstone2020faster}. Each experiment was run for
100 iterations; the largest distance over iterations
20 through 100 is reported.

\begin{table}[t]
\centering
\caption{Exact confinement examples with step size $\alpha=0.2$.
The improved radius is attained in each row. The comparison radius
$R_{\mathrm{JM}}$ follows from Johnstone and Moulin
\cite{johnstone2020faster} with weak sharp minima.}
\label{tab:exact}
\small
\begin{tabular}{lrrrr}
\toprule
Objective & $\mu$ & Observed & $\Rsharpbound$ & $R_{\mathrm{JM}}$\\
\midrule
$|x|$ & 1.0000 & 0.2000 & 0.2000 & 0.2236\\
$\|x\|_1$, $n=2$ & 1.0000 & 0.2828 & 0.2828 & 0.3464\\
$\|x\|_\infty$, $n=4$ & 0.5000 & 0.2000 & 0.2000 & 0.2828\\
$\|x\|_\infty$, $n=9$ & 0.3333 & 0.3000 & 0.3000 & 0.3606\\
$\|x\|_\infty$, $n=16$ & 0.2500 & 0.4000 & 0.4000 & 0.4472\\
$\|x\|_\infty$, $n=25$ & 0.2000 & 0.5000 & 0.5000 & 0.5385\\
\bottomrule
\end{tabular}
\end{table}

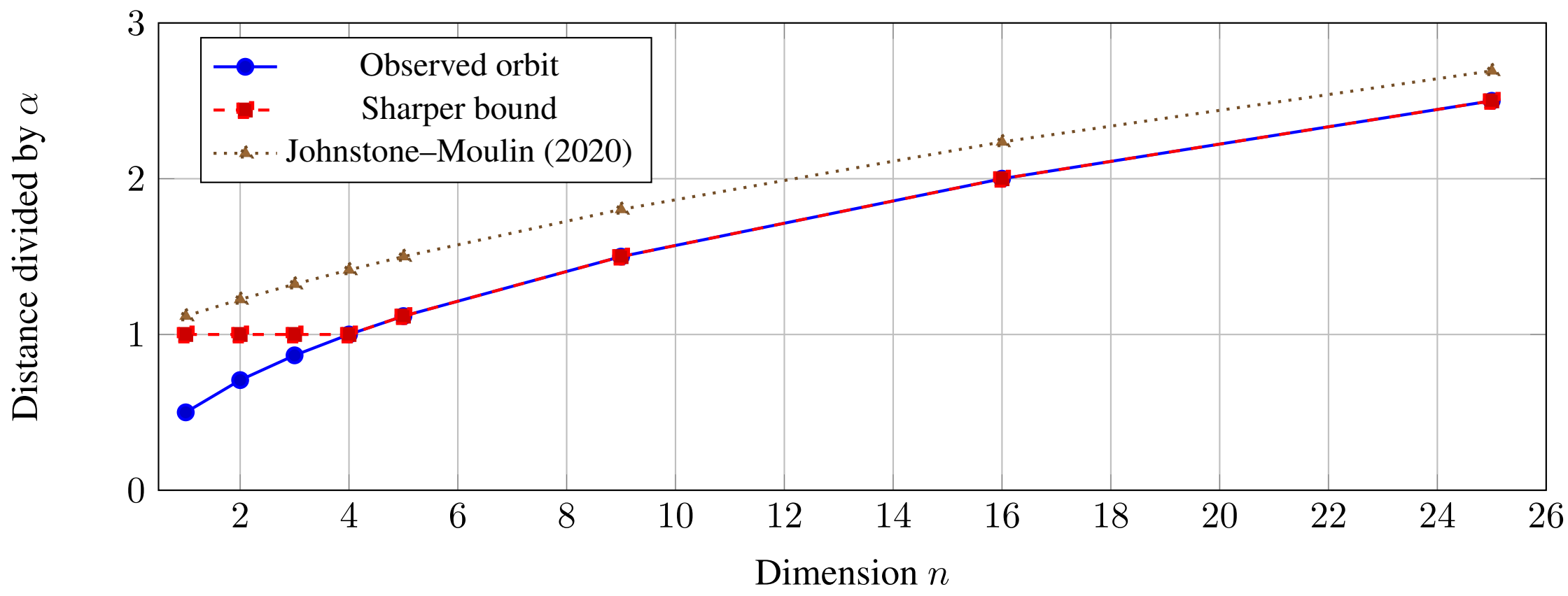
\begin{figure}[t]
\centering
\begin{tikzpicture}
\begin{axis}[
 width=0.89\linewidth,height=6.0cm,
 xlabel={Dimension $n$},ylabel={Distance divided by $\alpha$},
 xmin=0.5,xmax=26,ymin=0,ymax=3,
 grid=major,legend pos=north west,
 legend style={font=\small},
 every axis plot/.append style={thick}
]
\addplot+[mark=*,solid] coordinates {
 (1,0.5) (2,0.707107) (3,0.866025) (4,1) (5,1.118034)
 (9,1.5) (16,2) (25,2.5)};
\addlegendentry{Observed orbit}
\addplot+[mark=square*,dashed] coordinates {
 (1,1) (2,1) (3,1) (4,1) (5,1.118034)
 (9,1.5) (16,2) (25,2.5)};
\addlegendentry{Sharper bound}
\addplot+[mark=triangle*,dotted] coordinates {
 (1,1.118034) (2,1.224745) (3,1.322876) (4,1.414214)
 (5,1.5) (9,1.802776) (16,2.236068) (25,2.692582)};
\addlegendentry{Johnstone--Moulin (2020)}
\end{axis}
\end{tikzpicture}
\caption{Normalized confinement radii for $f(x)=\|x\|_\infty$
with $\alpha=0.2$. The Johnstone--Moulin radius is the weak-sharp
specialization of Theorem~1 in \cite{johnstone2020faster}.
The improved bound is attained for dimensions $n\geq4$.}
\label{fig:dimension}
\end{figure}

For $f(x)=\|x\|_\infty$, the Hausdorff distance is known
exactly: $H_\delta=\sqrt{n}\,\delta$. The published radius
of \cite{johnstone2020faster} is
$R_{\mathrm{JM}}=\alpha\sqrt{1+n/4}$, whereas the improved
radius is $\alpha\max\{1,\sqrt{n}/2\}$. Figure~\ref{fig:dimension}
shows both predictions together with the actual trajectory. For
$n\geq4$, the measured radius and improved radius agree exactly.
We repeated the experiments with step sizes $0.025$, $0.05$,
$0.1$, $0.2$, and $0.4$. The observed radius and both theoretical
radii scaled linearly with $\alpha$, as predicted.

\subsection{Anisotropic polyhedral objectives}
We next studied $f_\theta(x)=\max\{|a_1^Tx|,|a_2^Tx|\}$ where $a_1=(1,0)$ and $a_2=(\cos\theta,\sin\theta)$.
For $0<\theta\leq\pi/2$, the unique
minimizer is the origin, and $G=1$.
To determine the weak sharpness modulus,
let $u$ be the unit vector bisecting the
angle between $a_1$ and $a_2$, and let
$v$ be a perpendicular unit vector.
Every $x\in\R^2$ can be expressed as
$x=su+tv$, with $\norm{x}=\sqrt{s^2+t^2}$.
In these coordinates, $f_\theta(x) =\max\left\{ \left|s\cos(\theta/2)-t\sin(\theta/2)\right|, \left|s\cos(\theta/2)+t\sin(\theta/2)\right| \right\}$.
Using $\max\{|a-b|,|a+b|\}=|a|+|b|$, we obtain $f_\theta(x) = |s|\cos(\theta/2)+|t|\sin(\theta/2) \geq\sin(\theta/2)(|s|+|t|) \geq\sin(\theta/2)\norm{x}$.
Equality holds when $x=v$, so $\mu=\sin(\theta/2)$ is a weak sharpness modulus.

The function $f_\theta$ is positively homogeneous of degree one, meaning that $f_\theta(tx)=tf_\theta(x)$ for every $x\in\R^2$ and $t\geq0$.
We have already established that $f_\theta(x)\geq\mu\norm{x}$, so every $x\in S_\delta$ satisfies $\norm{x}\leq\delta/\mu$. Hence $H_\delta\leq\delta/\mu$. Furthermore, we identified a unit vector $v$ satisfying $f_\theta(v)=\mu$. By positive homogeneity, the point $x=(\delta/\mu)v$ satisfies $f_\theta(x)=\delta$ and $\norm{x}=\delta/\mu$. Therefore, $H_\delta=\delta/\mu$.

We considered four values of $\theta$: $0.20$, $0.35$, $0.60$, and $0.90$ radians. For each value, we generated 24 random initial points by drawing each coordinate independently from the uniform distribution on $[-3,3]$. We ran the subgradient method for 4000 iterations from each initial point with $\alpha=0.2$ and measured the largest distance over iterations 2000 through 4000. Table~\ref{tab:oblique} reports the largest late-iterate distance among 24 starts, together with the improved and Johnstone--Moulin radii. All observed distances were below both bounds.
The sharper bound is substantially more conservative than the actual trajectory for small $\theta$.

\begin{table}[t]
\centering
\caption{Oblique polyhedral objectives with $\alpha=0.2$,
24 starts per angle, and 4000 iterations per trajectory. The
comparison radius $R_{\mathrm{JM}}$ follows from
\cite{johnstone2020faster}. ``Observed'' is the maximum
late-iterate distance across the 24 starts.}
\label{tab:oblique}
\small
\begin{tabular}{rrrrr}
\toprule
$\theta$ (radians) & $\mu$ & Observed & $\Rsharpbound$ & $R_{\mathrm{JM}}$\\
\midrule
0.20 & 0.0998 & 0.1980 & 1.0017 & 1.0214\\
0.35 & 0.1741 & 0.1930 & 0.5744 & 0.6082\\
0.60 & 0.2955 & 0.1877 & 0.3384 & 0.3931\\
0.90 & 0.4350 & 0.1865 & 0.2299 & 0.3047\\
\bottomrule
\end{tabular}
\end{table}

\begin{figure}[t]
\centering
\begin{tikzpicture}
\begin{axis}[
 width=0.86\linewidth,height=5.9cm,
 xlabel={Angle $\theta$ (radians)},ylabel={Distance to minimizer},
 xmin=0.14,xmax=0.97,ymin=0,ymax=1.12,
 grid=major,legend pos=north east,
 legend style={font=\small},every axis plot/.append style={thick}
]
\addplot+[mark=*,solid] coordinates {
 (0.2,0.197998) (0.35,0.192974) (0.6,0.187674) (0.9,0.186483)};
\addlegendentry{Largest observed tail distance}
\addplot+[mark=square*,dashed] coordinates {
 (0.2,1.001669) (0.35,0.574356) (0.6,0.338386) (0.9,0.229903)};
\addlegendentry{Sharper bound}
\addplot+[mark=triangle*,dotted] coordinates {
 (0.2,1.021440) (0.35,0.608181) (0.6,0.393072) (0.9,0.304722)};
\addlegendentry{Johnstone--Moulin (2020)}
\end{axis}
\end{tikzpicture}
\caption{Measured distances and confinement radii for the
polyhedral objective. The
Johnstone--Moulin radius is implied by \cite{johnstone2020faster};
neither bound need be attained under the selected subgradients.}
\label{fig:oblique}
\end{figure}
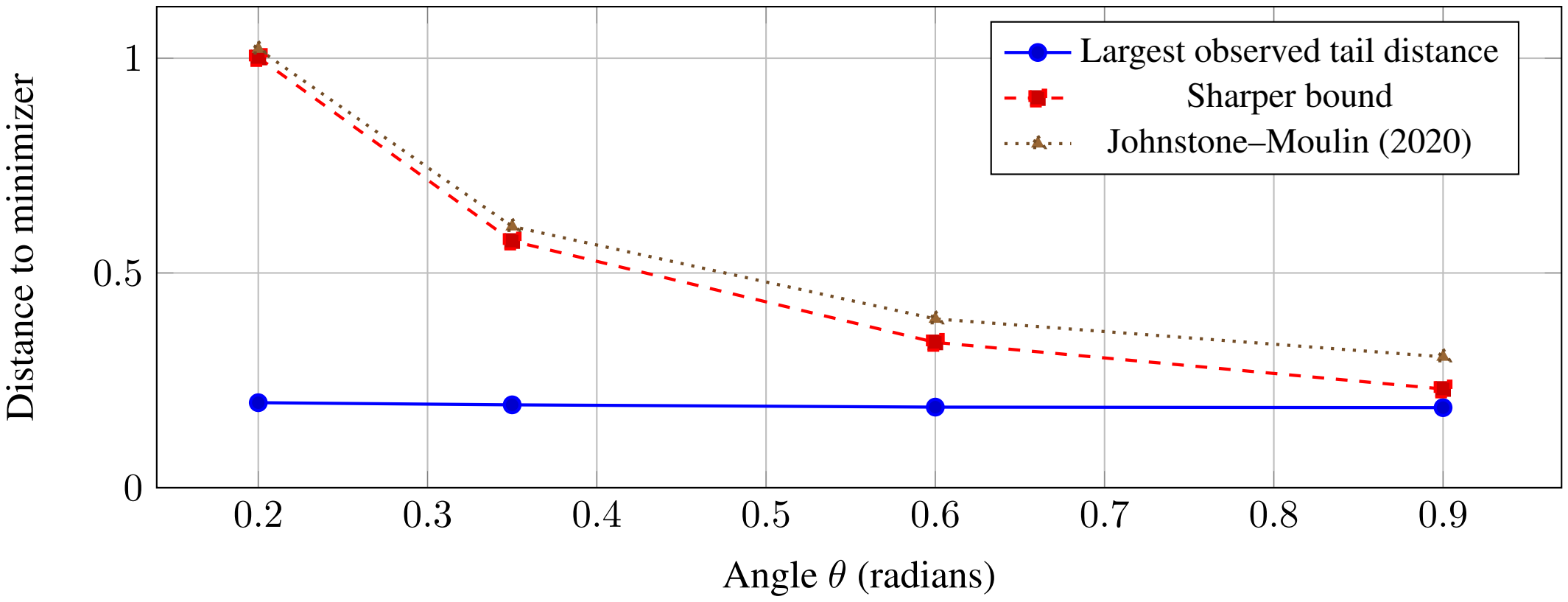

\subsection{Entrance time and an objective without a sharp minimum}
Corollary~\ref{cor:entranceattained} bounds the time required to first enter a near-optimal sublevel set when a minimizer exists.
For the polyhedral objective, we used $\theta=0.35$, $x_0=(2,-1)$, $\alpha=0.2$, and $\delta=0.15$. Here $H_\delta=0.8615$, and the confinement bound from Theorem~\ref{thm:hausdorff} is $\sqrt{H_\delta^2+\alpha^2G^2}=0.8844$.
Corollary~\ref{cor:entranceattained} guarantees entrance into $S_\delta$ in at most 250 iterations. The actual first entrance occurred at iteration 19, and the largest distance thereafter was $0.6717$ over a run of 4000 iterations. The finite-time entrance bound is a worst-case guarantee rather than a prediction of the actual entrance time.

Finally, we considered $f(x)=\sqrt{1+x^2}-1$, which is convex, continuously differentiable, and has $G=1$. Its unique minimizer is $0$ and its sublevel geometry is known exactly: $H_\delta=\sqrt{\delta(\delta+2)}$. This objective does not satisfy a global H\"olderian growth
condition of any exponent $\vartheta\in(0,1]$. Indeed, $f(x)\sim x^2/2$ as $x\to0$ and $f(x)\sim|x|$ as
$|x|\to\infty$. Thus an inequality $f(x)\geq c|x|^{1/\vartheta}$ cannot hold globally for any $c>0$ and $\vartheta\in(0,1]$. Consequently, the Johnstone--Moulin bounds do not apply to this example, while Theorem~\ref{thm:hausdorff} does.

For $0<\alpha\leq1$, let $\delta_0=\alpha/2$. Using the known expression for $H_\delta$ and the definition of $Q_{\delta,\alpha}$ in \eqref{eq:envelope}, the two asymptotic bounds from Corollary~\ref{cor:limit} are $\max\{H_{\delta_0},\sqrt{Q_{\delta_0,\alpha}}\} = \sqrt{\alpha+\alpha^2/4}$ and $\sqrt{H_{\delta_0}^2+\alpha^2} = \sqrt{\alpha+5\alpha^2/4}$.
The first bound retains the objective-error term through $Q_{\delta_0,\alpha}$, whereas the second uses the simpler upper bound in \eqref{eq:envelopeupper}. We ran 1200 iterations from $x_0=2$ at four step sizes. Table~\ref{tab:noheb} reports both bounds and the largest distance over iterations 1000 through 1200.

\begin{table}[t]
\centering
\caption{Comparison of the observed late-iterate distances with the two theoretical confinement bounds for $f(x)=\sqrt{1+x^2}-1$. The observed distance is the maximum over iterations 1000 through 1200.}
\label{tab:noheb}
\begin{tabular}{rrrr}
\toprule
$\alpha$ & $\max_{1000\leq k\leq1200}|x_k|$ & $\max\{H_{\delta_0},\sqrt{Q_{\delta_0,\alpha}}\}$ & $\sqrt{H_{\delta_0}^2+\alpha^2}$\\
\midrule
0.05 & $<10^{-16}$ & 0.2250 & 0.2305\\
0.10 & $<10^{-16}$ & 0.3202 & 0.3354\\
0.20 & $<10^{-16}$ & 0.4583 & 0.5000\\
0.40 & $<10^{-16}$ & 0.6633 & 0.7746\\
\bottomrule
\end{tabular}
\end{table}
The iterates converge in this smooth example; the bounds are conservative but remain valid.

\section{Discussion}
Our results distinguish three levels of assumptions for a constant-step projected subgradient trajectory. First, if a continuous convex objective has a finite infimum on a nonempty, closed, convex feasible set and the selected subgradients are bounded, the trajectory enters a near-optimal sublevel set in finite time. This conclusion does not require a minimizer, a growth condition, or any Hausdorff-distance assumption. The entrance-time estimate uses an approximately optimal feasible point; when a minimizer exists, comparison with the optimal set provides a simpler, stronger bound.

Second, if the infimum is attained, distance to the optimal set is well defined. Theorem~\ref{thm:hausdorff} bounds the distance of every iterate after the first entrance by a function of the Hausdorff distance between the near-optimal sublevel set and the optimal set. This inequality requires no H\"olderian growth condition and remains formally valid if the Hausdorff distance is infinite. Finiteness of that distance is what makes the resulting confinement radius finite. Retaining the objective-error term in the one-step subgradient inequality can improve the radius further. A compact optimal set is sufficient, but not necessary, for the finite Hausdorff-distance property.

Third, the additional assumption of weak sharp minima yields the explicit asymptotic distance estimate of Theorem~\ref{thm:sharp}. This improves the radius implied by Johnstone and Moulin \cite{johnstone2020faster}. Section~\ref{sec:attainmentExamples} constructs trajectories that attain the bound, while the oblique experiments show that a worst-case bound may be conservative for a particular objective and subgradient selection rule. The smooth example demonstrates that a finite Hausdorff confinement radius can be obtained without a global H\"olderian growth condition.

The analysis is deterministic and applies to the projected subgradient method, which includes the subgradient method when $C=\R^n$. Bounded selected subgradients are assumed for all three results; the existence of minimizers and weak sharpness are needed only where explicitly stated.

In future work, we hope to create an adaptive step size rule for subgradient methods that takes advantage of the selected subgradient magnitudes to achieve fast convergence.


\end{document}